\documentclass[twoside,draft,reqno]{amsart}
\usepackage{amssymb,amsfonts,amsxtra}
\usepackage[all]{xy}
\usepackage{fullpage}
\newtheorem{theorem}{Theorem}[section]
\newtheorem{cor}[theorem]{Corollary}
\newtheorem{lem}[theorem]{Lemma}
\newtheorem{prop}[theorem]{Proposition}
\newtheorem{defi}[theorem]{Definition}
\newtheorem{rem}[theorem]{Remark}
\newtheorem{con}[theorem]{Conjecture}
\begin{document}
\title{Spaces of multiplicative maps between highly structured ring spectra}
\author{A. Lazarev}
\address{Mathematics Department, Bristol University, Bristol, BS8 1TW,
England.} \email{a.lazarev@bristol.ac.uk}
\begin{abstract}
We uncover a somewhat surprising connection between spaces of
multiplicative maps between $A_\infty$-ring spectra and
topological Hochschild cohomology. As a consequence we show that
such spaces become infinite loop spaces after looping only once.
We also prove that any multiplicative cohomology operation in
complex cobordisms theory $MU$ canonically lifts to an
$A_\infty$-map $MU\rightarrow MU$.  This implies, in particular,
that the Brown-Peterson spectrum $BP$ splits off $MU$ as an
$A_\infty$-ring spectrum.
\end{abstract}
\maketitle
\section{Introduction}  The main purpose of the
present work is to provide a workable method for computing the
homotopy type of spaces of $A_\infty$-maps between $A_\infty$-ring
spectra (or $S$-algebras in the terminology of \cite{EKMM}). We
make substantial use of the previous results by the author in
\cite{Laz} and for the reader's convenience a brief summary of
these is given in Section $2$. In Section $3$ we collect
miscellaneous technical results concerning function spectra and
topological Hochschild cohomology. Some of these results are
surely known to experts but never have been written down. The
formula (\ref{mod}) (base change) deserves special mention. While
easy to prove it is extremely convenient  when computing with
various spectral sequences.

Our first main result (Theorem \ref{osn}) essentially states that
the mapping space between two $S$-algebras $A$ and $B$ is
determined after taking based loops by the spectrum of topological
derivations ${\bf Der}(A,B)$. Therefore the problem of computing
higher homotopy groups of this mapping space is a problem of
stable homotopy which turns out to be quite amenable, particularly
because in many cases ${\bf Der}(A,B)$ can be reduced to ${\bf
THH}(A,B)$, the topological Hochschild cohomology of $A$ with
values in $B$.

The computation of the zeroth homotopy group of the mapping space
is, of course, a completely different story. We give a simple
answer in the special case when $A$ is a connective $S$-algebra
while $B$ is coconnective (that is, with vanishing homotopy in
positive dimensions). This is Theorem \ref{eka}.

Even though the problem of computing homotopy classes of
$S$-algebra maps $A\rightarrow B$ is essentially unstable it does
lend itself to analysis by methods of obstruction theory developed
in \cite{Laz}. Our second main result (Theorem \ref{maps})
demonstrates that any multiplicative cohomology operation in
complex cobordism theory $MU$ canonically (even uniquely in an
appropriate sense) lifts to an $S$-algebra self-map of $MU$. This
is used to show that for an $S$-algebra $E$ belonging to a fairly
large class of complex-oriented theories any multiplicative
operation $MU\rightarrow E$ lifts to an $S$-algebra map. Another
corollary is that the Brown-Peterson spectrum $BP$ splits off $MU$
localized at $p$ as an $S$-algebra.

The paper is written in the language of $S$-modules of \cite
{EKMM} and we routinely use the results and terminology of the
cited reference.

{$ Notations.$} In Sections $2$ and $3$  we work in the category
of modules or algebras over a fixed $q$-cofibrant commutative
$S$-algebra $R$, the smash product $\wedge$ and the function
spectrum $F(-,-)$ are always understood as $\wedge_R$ and
$F_R(-,-)$. In Section $4$ we specialize to $R=S$. The free
$R$-algebra on an $R$-module $M$ is denoted by $T(M)$. The space
of maps between two $R$-algebras $A$ and $B$ is denoted by
$F_{R-alg}(A,B)$.  If $A$ and $B$ are commutative $R$-algebras
then $[A,B]_c$ denotes the set of homotopy classes of commutative
algebra maps from $A$ to $B$. If $A$ and $B$ are associative
$R$-algebras then $[A,B]_a$ stands for homotopy classes of
associative algebra maps. For an associative $R$-algebra $A$ and
two $A$-modules $M$ and $N$ we denote by $[M,N]_{A-mod}$ homotopy
classes of $A$-module maps from $M$ to $N$. Similarly
$[M,N]_{A-bimod}$ stands for homotopy classes of bimodule maps.
Finally we denote by $Mult(E,F)$ the set of multiplicative (up to
homotopy) maps between ring spectra $E$ and $F$.
 To distinguish between strict
isomorphisms and weak equivalences we will use, as a rule, the
symbol $`\cong$' for the former and $`\simeq$' for the latter.

\section{Topological derivations and topological singular extension of $S$-algebras}
In this short section we give an overview of some of the author's
results from \cite{Laz} which will be needed later on.

Let $A$ be a $q$-cofibrant  $R$-algebra and $M$ a $q$-cofibrant
$A$-bimodule. Then the $R$-module $A\vee M$ has the obvious
structure of an $R$-algebra (`square-zero extension' of $A$).
Consider the set $[A,A\vee M]_{a/A}$ of homotopy classes of
$R$-algebra maps from $A$ to $A\vee M$ in the category of
$R$-algebras over $A$, that is the $R$-algebras supplied with an
$R$-algebra map into A.
\begin{theorem}\label{lp}
There exists an $A$-bimodule $\Omega_A$ and a natural in $M$
isomorphism $$[A,A\vee M]_{a/A}\cong [\Omega_A,M]_{A-bimod}$$
 where the right hand side denotes the homotopy classes of maps in the category
of $A$-bimodules.
\end{theorem}
\begin{rem}Sometimes we will need a refinement of the above
theorem which is formulated as follows. Let $B$ be an $R$-algebra
over $A$, i.e. there exists a fixed $R$-algebra map $B\rightarrow
A$. Then there is a natural isomorphism $$[A,B\vee M]_{a/B}\cong
[A\wedge_B\Omega_B\wedge_BA,M]_{A-bimod}.$$ Furthermore an
$A$-bimodule $M$ can be considered as a $B$-bimodule and we have
$$[B,B\vee M]_{a/B}\cong [\Omega_B,M]_{B-bimod}\cong
[A\wedge_B\Omega_B\wedge_BA,M]_{A-bimod}\cong[B,A\vee M]_{a/A}.$$
The isomorphism $[B,B\vee M]_{a/B}\cong [B,A\vee M]_{a/A}$ will be
used without explicit mention later on in this paper.\end{rem}
\begin{defi}
The topological derivations  $R$-module of $A$ with values in $M$
is the function $R$-module $F_{A\wedge A^{op}}(\Omega_A,M)$. We
denote it by ${\bf Der}_R(A,M)$ and its $i$th homotopy group by
$Der^{-i}_R(A,M)$.
\end{defi}
The $A$-bimodule $\Omega_A$ is constructed as the $q$-cofibrant
approximation of the homotopy fibre of the multiplication map
$A\wedge A\rightarrow A$. There exists the following homotopy
fibre sequence of $R$-modules:
\begin{equation}\label{sdf}{\bf THH}_R(A,M)\rightarrow M \rightarrow
{\bf Der}_R(A,M)\end{equation} Here ${\bf THH}_R(A,M)$ is the
topological Hochschild cohomology spectrum of $A$ with values in
$M$: \[{\bf THH}_R(A,M):=F_{A\wedge A^{op}}(\tilde{A},M)\] where
$\tilde{A}$ is the $q$-cofibrant replacement of $A$ as an
$A$-bimodule.

We will also have a chance to use topological Hochschild {\it
homology} spectrum ${\bf THH}^R(A,M):=A\wedge_{A\wedge A^{op}}M$.
If the $R$-algebra $A$ is commutative and  the left and right
$A$-module structures on $M$ agree
 then both ${\bf THH}^R(A,M)$
and ${\bf THH}_R(A,M)$ are $A$-modules and there is a weak
equivalence of $A$-modules \[{\bf THH}_R(A,M)\cong F_A({\bf
THH}^R(A,M),A).\] Furthermore in this case the sequence
(\ref{sdf}) splits giving a canonical weak equivalence ${\bf
THH}_R(A,M)\simeq \Sigma^{-1}{\bf Der}_R(A,M)\vee M.$

Suppose we are given a topological derivation $d:A\rightarrow
A\vee M$. Consider the following homotopy pullback diagram of
$R$-algebras \[\xymatrix{X\ar[r]\ar[d]&A\ar[d]\\ A\ar^-d[r]&A\vee
M }\] Here the rightmost downward arrow is the canonical inclusion
of a retract. Then we have the following homotopy fibre sequence
of $R$-modules:
\begin{equation}\label{er}\xymatrix{\Sigma^{-1}M\ar[r]& X\ar[r]& A}\end{equation}
\begin{defi} The homotopy fibre sequence (\ref{er}) is called the topological
singular extension associated with the derivation  $d:A\rightarrow
A\vee M$.
\end{defi}
\begin{theorem}\label{ext} Let $\Sigma^{-1}M\rightarrow X\rightarrow A$ be a singular
extension of $R$-algebras associated with a derivation
$d:A\rightarrow A\vee M$ and $f:B\rightarrow A$ a map of
$R$-algebras. Then $f$ lifts to an $R$-algebra map $B\rightarrow
X$ iff a certain element in $Der^0_R(B,M)$ is zero. Assuming that
a lifting exists  the homotopy fibre of the map
\[F_{R-alg}(B,X)\rightarrow F_{R-alg}(B,A)\]
over the point $f\in F_{R-alg}(B,A)$ is weakly equivalent to
$\Omega^\infty {\bf Der}_R(B,\Sigma^{-1}M)$, the $0$th space of
the spectrum ${\bf Der}_R(B,\Sigma^{-1}M)$.
\end{theorem}
\begin{theorem}\label{pos} Assume that $R$ is  connective and $A$
is a connective $R$-algebra. Then the Postnikov tower of $A$
\[\xymatrix{A_0=H\pi_0A&A_1\ar[l]&\ldots\ar[l]&A_n\ar[l]&
A_{n+1}\ar[l]& \ldots\ar[l]}\] is a tower of $R$-algebras.
Moreover the homotopy fibre sequences \[\xymatrix{ A_n&
A_{n+1}\ar[l]&\Sigma^{n+1}H\pi_{n+1}A\ar[l]}\] are topological
singular extensions.
\end{theorem}

\section{Base change and topological derivations of supplemented $R$-algebras }
In this section we discuss topological Hochschild cohomology and
topological derivatons of supplemented $R$-algebras and the
behaviour of the forgetful map $l$ in the hypercohomology spectral
sequence. This material is largely parallel to \cite{Laz}, section
$9$ and so most of the proofs will be omitted.

We'll start with some general lemmas.
\begin{lem}\label{ol}
Let $A$, $B$, $C$ be $R$-algebras and $M$, $N$, $L$ be an $A\wedge
B$-module, an $C\wedge A^{op}$-module and a $C\wedge B$-module
respectively. Then there is a natural isomorphism of $R$-modules:
$$F_{A\wedge B}(M, F_{C}(N,L))\cong F_{C\wedge B}(N\wedge_A
M,L).$$
\end{lem}
\begin{proof} Let us first check the above equivalences for $M=A\wedge
B\wedge \tilde{M}$ and $N=C\wedge A^{op}\wedge\tilde{N}$. We have:
\begin{align*}F_{A\wedge B}(M,F_C(N,L))&\cong F_{A\wedge B}(A\wedge
B\wedge\tilde{M}, F_C(C\wedge A^{op}\wedge \tilde{N},L))\\&\cong
F(\tilde{M},F(A^{op}\wedge\tilde{N},L))\\& \cong F(\tilde{M}\wedge
A^{op}\wedge\tilde{N},L).\end{align*} Likewise,
\begin{align*}F_{C\wedge B}(N\wedge_A M, L)&\cong F_{ C\wedge B}(C\wedge
A^{op}\wedge\tilde{N} \wedge_A A\wedge B\wedge \tilde{M},
L)\\&\cong F(A^{op}\wedge\tilde{N}\wedge\tilde{M}, L).\end{align*}
Observe that the above isomorphisms are natural in $M$ and $N$
that is, with respect to arbitrary maps of $A\wedge B$ modules
$A\wedge B\wedge \tilde{M}_1\rightarrow A\wedge B\wedge
\tilde{M}_2$ and of $C\wedge A^{op}$-modules $C\wedge A^{op}\wedge
\tilde{N}_1\rightarrow C\wedge A^{op}\wedge \tilde{N}_2$ (not only
those coming from $\tilde{M}_1\rightarrow \tilde{M}_2$ and
$\tilde{N}_1\rightarrow \tilde{N}_2$). To obtain the general case
it suffices to notice that for any $M$ and $N$ there exist
standard split coequalizers of $R$-modules $$\xymatrix{(A\wedge
B)^{\wedge 2}\wedge M \ar@<-0.5ex>[r]\ar@<0.5ex>[r]&A\wedge
B\wedge M\ar[r]&M }$$ and $$\xymatrix{(C\wedge A^{op})^{\wedge
2}\wedge N \ar@<-0.5ex>[r]\ar@<0.5ex>[r]&C\wedge A^{op}\wedge N
\ar[r]&N }.$$ With this Lemma \ref{ol} is proved. \end{proof}

Now  let $C=N$ be an $A$-bimodule via an $R$-algebra map
$f:A\rightarrow C$. Then the $C\wedge B$-module $L$ acquires a
structure of an $A\wedge B$-module via the map \[\xymatrix{A\wedge
B\ar^{f\wedge id}[r]& C\wedge B}.\] Furthermore simple diagram
chase shows that the isomorphism $F_C(C,L)\cong L$ is in fact an
isomorphism of $A\wedge B$-modules. This gives the following
\begin{cor}\label{kar}  There exists the following natural
isomorphism of $R$-modules:
\begin{equation}\label{mod}
F_{C\wedge B}(C\wedge_A M,L)\cong F_{A\wedge B}(M,F_{C}(C,L))\cong
F_{A\wedge B}(M,L)
\end{equation}
\end{cor}
 We will refer to the
isomorphism (\ref{mod}) as base change. Related formulae are found
in \cite{EKMM}, III.6.

\begin{cor}\label{alg} If an $R$-algebra $B$ is an $A$-bimodule via a
$q$-cofibration of $R$-algebras $A\rightarrow B$, then ${\bf
THH}_R(A,B)\simeq F_{{B}\wedge A^{op}}(\tilde{B},\tilde{B})$ where
$\tilde{B}$ is the $q$-cofibrant approximation of the $B\wedge
A^{op}$-module $B$. In particular, ${\bf THH}_R(A,B)$ is an
$R$-algebra under the composition product.\end{cor} \begin{proof}
Denoting by $\tilde{A}$ the $q$-cofibrant approximation of the
$A$-bimodule $A$ we have the following isomorphisms of
$R$-modules:
\begin{align*}{\bf THH}_R(A,B)&\cong F_{A\wedge
A^{op}}(\tilde{A},B)\\ &\cong F_{\tilde{A}\wedge A^{op}}
(A,F_B(B,B))\\ &\cong F_{B\wedge A^{op}}(B\wedge_A
\tilde{A},B).\end{align*} The $B\wedge A^{op}$-module $B\wedge_A
\tilde{A}\cong B\wedge A^{op}\wedge_{A\wedge A^{op}}\tilde{A}$ is
a $q$-cofibrant $B\wedge A^{op}$-module because the functor
$?\rightarrow B\wedge A^{op}\wedge_{A\wedge A^{op}}?$ preserves
$q$-cofibrant modules. Therefore $ F_{B\wedge A^{op}}(B\wedge_A
\tilde{A},B)$ represents derived function $B\wedge A^{op}$-module
and is equivalent to $F_{{B}\wedge A^{op}}(\tilde{B},\tilde{B})$.
\end{proof}
We now discuss topological Hochschild cohomology and derivations
of supplemented $R$-algebras. Let $A$ be a $q$-cofibrant
$R$-algebra. We say
 that $A$ is supplemented if it is supplied with an $R$-algebra morphism
$\epsilon:A\rightarrow B$ which we will assume to be a
$q$-cofibration of $R$-algebras. Denote by $\Omega^B_A$ the
homotopy fibre of the map $B\wedge A\rightarrow B$ that determines
the structure of a right $A$-module on $B$. We will assume without
loss of generality that $\Omega^B_A$ is a $q$-cofibrant right
$A$-module. Recall that the module of differentials $\Omega_A$ for
$A$ is defined from the homotopy fibre sequence
$$\Omega_A\rightarrow A\wedge A\rightarrow A$$ where the second
arrow is the multiplication map. Smashing this fibre sequence on
the left with $B$ over $A$ we get the fibre sequence $$
B\wedge_A\Omega_A\rightarrow B\wedge A\rightarrow B$$ That shows
that $\Omega^B_A$ is weakly equivalent to $B\wedge_A\Omega_A$ as a
$B\wedge A$-module. Further, base change gives a weak equivalence
$$F_{B\wedge A^{op}}(\Omega^B_A,B)\simeq F_{B\wedge
A^{op}}(B\wedge_A\Omega_A,B)\simeq F_{A\wedge
A^{op}}(\Omega_A,B)\cong {\bf Der}_R(A,B).$$

Recall from \cite{Laz} that there is a `universal derivation'
$d:A\rightarrow \Omega_A$ which is defined as the composite map
$$A\rightarrow A\vee \Omega_A\rightarrow \Omega_A$$ where the
first map is the map of algebras over $A$ adjoint to the identity
map $\Omega_A\rightarrow \Omega_A$ and the second map is the
projection. The universal derivation allows one to define the
forgetful map $l:{\bf Der}_R(A,B)\rightarrow F(A,B)$ as the
composite map $${\bf Der}_R(A,B)\simeq F_{A\wedge
A^{op}}(\Omega_A,B)\rightarrow F(\Omega_A,B) \rightarrow F(A,B)$$
where the last map is induced by $d$. In terms of $\Omega^B_A$ the
forgetful map $l$ admits the following description. The fibre
sequence $\Omega^B_A\rightarrow B\wedge A\rightarrow B$ splits via
the map $B\cong B\wedge R\stackrel{id\wedge
1}{\longrightarrow}B\wedge A$ so there is a weak equivalence of
$R$-modules $B\wedge A\simeq B\vee \Omega^B_A$.
 Denote by $\overline{d}:A\rightarrow \Omega^B_A$ the following
 composite map
 \[\xymatrix{A\ar^-{\epsilon\wedge id}[r]&B\wedge A\simeq
 B\vee \Omega^B_A\ar[r]&\Omega^B_A}\]
 the last arrow being the projection onto the wedge summand.
Then $l$ coincides with the following composition:
\[\xymatrix{{\bf Der}_R(A,B)=F_{B\wedge A^{op}}(\Omega^B_A,B)
\ar[r]& F(\Omega^B_A,B)\ar[r]& F(A,B)}\] the first arrow being the
forgetful map and the second one is induced by $\overline{d}$.

Next we discuss the behaviour of the forgetful map $l$ in the
hypercohomology spectral sequence. To do this we need to review
algebraic Hochschild cohomology for supplemented algebras. The
exposition will be somewhat sketchy since it is parallel to
\cite{Laz}, section $9$.
\begin{defi} Let $A_\ast$ be a graded algebra over a graded commutative
algebra $R_\ast$ supplied with an $R_\ast$-algebra map
 $\epsilon:A_\ast\rightarrow B_\ast$ (supplementation). Then algebraic
Hochschild cohomology of $A_\ast$ with coefficients in $B_\ast$ is
defined as $$HH^{\ast}_{R_\ast}(A_\ast, B_\ast)=
Ext^\ast_{B_\ast\otimes_{R_\ast}^LA_\ast^{op}}(B_\ast, B_\ast)$$
where $\otimes_{R_\ast}^L$ denotes the derived tensor product
\end{defi}
\begin{rem} If $A_\ast$ is flat as an $R_\ast$-module, then this definition is
equivalent to the standard one found in, e.g. \cite{CaE}
\end{rem}
 We also have a generalization of the standard complex
which computes Hochschild cohomology. Let $\tilde{A}_\ast$ be a
differential graded supplemented $R_\ast$-algebra which is
quasiisomorphic to ${A}_\ast$ and $R_\ast$-projective. Denote by
$\tilde{\epsilon}:\tilde{A}_\ast\rightarrow B_\ast$ its
supplementation. Consider the bar-resolution of the right
$\tilde{A}_\ast$-module $B_\ast$; here and later on $\otimes$
stands for $\otimes_{R_\ast}$
\begin{equation}\label{bnm}\xymatrix{B_\ast&B_\ast\otimes\tilde{A}_\ast
\ar[l]& B_\ast\otimes\tilde{A}_\ast\otimes\tilde{A}_\ast
\ar[l]&\ldots\ar[l]}
\end{equation}
with the usual bar differential
\[\partial(b\otimes\tilde{a_1}\otimes\ldots\otimes\tilde{a_n})= \pm
b\tilde{\epsilon}(\tilde{a}_1) \otimes\tilde{a}_2\ldots
\otimes\tilde{a}_n+ \Sigma \pm b\otimes\tilde{a}_1\otimes\ldots
\otimes \tilde{a}_i\tilde{a}_{i+1}\otimes\ldots\otimes\tilde{a}_n
.\] We don't specify the signs in this well-known formula, see
e.g. \cite{Mac}, Chapter X.
 This is actually a bicomplex since $\tilde{A}$ is a
differential graded algebra. Applying the functor
$Hom_{B_\ast\otimes\tilde{A}_\ast^{op}}(?,B_\ast)$ to (\ref{bnm})
we get the standard Hochschild cohomology (bi)complex
\[C^{ij}(A_\ast, B_\ast)=Hom^i(\tilde{A}_\ast^{\otimes
j},B_\ast).\] Now define the module of differentials
$\Omega_{A_\ast}^{B_\ast}$ from the following short exact sequence
\[\xymatrix{0\ar[r]& \Omega_{A_\ast}^{B_\ast}\ar[r]& B_\ast
\otimes\tilde{A}_\ast\ar[r]&B_\ast\ar[r]& 0}.\] Clearly
$\Omega_{A_\ast}^{B_\ast}$ is quasiisomorphic as a complex of
right $\tilde{A}_\ast$-modules to the truncated bar-resolution:
\begin{equation}\label{yui}\xymatrix{B_\ast\otimes\tilde{A}_\ast^{\otimes 2}
& B_\ast\otimes\tilde{A}_\ast^{\otimes 3}\ar[l]&\ldots\ar[l]}
\end{equation}
The universal derivation $d:\tilde{A}_\ast\rightarrow
\Omega_{A_\ast}^{B\ast}$ is induced by the map
$\tilde{A}_\ast\rightarrow B_\ast\otimes\tilde{A}_\ast$,
$\tilde{a}\rightarrow
1\otimes\tilde{a}-\tilde{\epsilon}(\tilde{a})\otimes 1$. If we
take the complex (\ref{yui}) as a model for
$\Omega_{A_\ast}^{B_\ast}$ then the universal derivation $d$ is a
map of complexes \[\xymatrix{d:\tilde{A}\ar[r]& \{
B_\ast\otimes\tilde{A}_\ast^{\otimes 2} \rightarrow
B_\ast\otimes\tilde{A}_\ast^{\otimes 3}\rightarrow\ldots\}}\]
where $\tilde{A}$ is considered to be a complex concentrated in
degree $0$ and $d(\tilde{a})=-1\otimes \tilde{a}\otimes 1$.

Further define algebraic derivations of $A_\ast$ with coefficients
in $B_\ast$ as $$Der^\ast_{R_\ast}(A_\ast, B_\ast)=
Ext^\ast_{B_\ast\otimes
\tilde{A}^{op}_\ast}(\Omega_{A_\ast}^{B_\ast}, B_\ast).$$ Then the
(truncated) standard resolution (\ref{yui}) provides a (bi)complex
for computing $Der^\ast_{R_\ast}(A_\ast, R_\ast)$:
\[\xymatrix{\overline{C}^\ast({A_\ast},B_\ast): Hom(\tilde{A_\ast},B_\ast)
\ar[r]& Hom(\tilde A_\ast^{\otimes 2},B_\ast) \ar[r]&\ldots.}\]
(This is  indeed a bicomplex, the additional differential being
induced from the internal differential in $\tilde{A_\ast}$). As in
the topological case the universal derivation determines the
forgetful map \[\xymatrix{l_{alg}: Der^\ast_{R_\ast}(A_\ast,
B_\ast)= Ext^\ast_{B_\ast\otimes \tilde{A}^{op}_\ast}
(\Omega_{A_\ast}^{B_\ast}, B_\ast)\ar[r]&
Hom^\ast(\tilde{A_\ast},B_\ast).}\]

Then we have the obvious
\begin{prop} The forgetful map
$$Der^\ast_{R_\ast}(A_\ast, B_\ast)\rightarrow
Hom^\ast(\tilde{A_\ast},B_\ast)=
Ext^\ast_{R_\ast}(A_\ast,B_\ast)$$ is induced by the projection
$\overline{C}^\ast({A_\ast},B_\ast)\rightarrow
Hom^\ast(\tilde{A_\ast},B_\ast)$ times $(-1)$
\end{prop}

Returning to our topological situation we have the following
result which is analogous to Proposition $9.3$ in \cite{Laz}:
\begin{prop}\label{an}
Let $A$, $B$ be $R$-algebras, $A\rightarrow B$ is an $R$-algebra
map. Suppose that the Kunneth spectral sequence for
$\pi_\ast(B\wedge A^{op})$ collapses and
 there is a ring isomorphism
$$\pi_{\ast}(B\wedge A^{op})\cong B_\ast\otimes^L_{R_\ast}
A^{op}_\ast$$ Then there are the following spectral sequences
$$^1E^{\ast\ast}_2=Der^\ast_{R_\ast}(A_\ast,B_\ast)=
Ext^\ast_{Tor^{R_\ast}_\ast(B_\ast,
A_\ast^{op})}(\Omega_{A_\ast}^{B_\ast},R_\ast) \Rightarrow {
Der}^\ast(A, B);$$
$$^2E^{\ast\ast}_2=Ext^\ast_{R_\ast}(A_\ast,B_\ast)\Rightarrow
[A,B]^\ast.$$ Furhermore, the forgetful map $l:
Der^\ast(A,B)\rightarrow [A,B]^\ast$ induces a map of spectral
sequences $^1E^{\ast\ast}_\ast\rightarrow~^2E^{\ast\ast}_\ast$
which on the level of $E_2$-terms gives the forgetful map
$l_{alg}:=Der^\ast_{R_\ast}(A_\ast,B_\ast)\rightarrow
Ext^\ast_{R_\ast}(A_\ast,B_\ast).$
\end{prop}
\section{Mappings spaces via derivations}
In this section we show that for two $R$-algebras $A$ and $B$ the
higher homotopy groups of the space $F_{R-alg}(A,B)$ can be
reduced to the computation of certain topological derivations.
This is important because in many cases topological derivations
can be further reduced to topological Hochschild cohomology which
is an essentially stable object, so that one could apply standard
methods of homological algebra for computation. As usual, we
assume that $A$ is a $q$-cofibrant $R$-algebra.

Now consider the $R$-module $A\vee\Sigma^{-d}A$, $d>0$. It can be
supplied canonically with the structure of an $R$-algebra over $A$
so that $\Sigma^{-d}A$ is a `square-zero ideal'. Denote this
$R$-algebra by $A_d$.

Let us also introduce the algebra $A(d):= A^{S^d}$, the cotensor
of A and the $d$-sphere $S^d$. Then $A(d) \cong F({\bf
R}\Sigma^{\infty}S^d_+,A)$ as an $R$-module (here ${\bf R}$ stands
for the free $R$-module functor. The structure of an $R$-algebra
on $A(d)$ is induced by the $R$-algebra structure on $A$ and the
topological diagonal $S^d\rightarrow S^d\times S^d$. The
coefficient rings of $A(d)$ and $A_d$ are both isomorphic to the
exterior algebra $\Lambda_{A_\ast}(y)$ where $y$ has degree $-d$.
There is also a weak equivalence of $R$-modules: \[A(d)\simeq
A\vee \Sigma^{-d}A\cong A_d.\] Notice that both $A(d)$ and $A_d$
are $R$-algebras over $A$, that is there exist maps of
$R$-algebras $A(d)\rightarrow A$ and $A_d\rightarrow A$. (The
first map is induced by choosing a base point in $S^d$, the second
map is the canonical projection).
\begin{theorem}\label{thy}The $R$-algebras $A(d)$ and $A_d$ are weakly equivalent
in the category of $R$-algebras.
\end{theorem}
\begin{proof} First consider the case $A=R$. Since $R$ is an
$R(d)$-module it makes sense to consider self-maps of $R$ in the
category of $R(d)$-modules. Notice that $R(d)$ is actually a
commutative $R$-algebra so we need not distinguish between left
and right $R(d)$-modules..
\begin{lem}\label{cop}
$$\pi_\ast F_{R(d)}(R,R)=R_\ast[[x]]$$ where the element x has
degree $d-1$.
\end{lem}
\begin{proof} Assume first that $R=S$, the sphere spectrum. We need
this special case because the connectiveness of $S$ will be used.
If $R$ is connective this step could be skipped. Consider the
spectral sequence
\begin{align*}E_{2}^{\ast\ast}&=Ext^{\ast\ast}_{S(d)_\ast}(S_\ast, S_\ast)\\&=
Ext^{\ast\ast}_{\Lambda_{S_\ast}(y)}(S_\ast,
S_\ast)\\&=S_\ast[[x]] \Longrightarrow
\pi_{\ast}F_{S(d)}(S,S).\end{align*} Here the element $x$ has
degree $d-1$ This spectral sequence collapses for dimensional
reasons. By Boardman's criterion \cite{Boa} it converges strongly
to its target which is complete with respect to the (cobar)
filtration. Since this filtration coincides with the $x$-adic
filtration on the associated graded $S_\ast$-module we conclude
that $$\pi_\ast F_{S(d)}(S,S)=S_\ast[[x]].$$ Notice that the fact
that $S$ is connective was used to show the collapse of our
spectral sequence. For instance if $d=1$ the elements $x^k$ are
located along the line of slope $1$ and the whole spectral
sequence $E_{2}^{\ast\ast}$
 lies above it.

Now let $R$ be an arbitrary commutative $S$-algebra. Consider the
spectral sequence
\begin{align*}^{\prime}E_{2}^{\ast\ast}&=Ext^{\ast\ast}_{R(d)_\ast}(R_\ast,
R_\ast)\\ &= Ext^{\ast\ast}_{\Lambda_{R_\ast}(y)}(R_\ast,
R_\ast)\\&=R_\ast[[x]] \Longrightarrow
\pi_{\ast}F_{R(d)}(R,R)\end{align*} and notice that the unit map
$S\rightarrow R$ determines the map of spectral sequences
$E_{2}^{\ast\ast}\rightarrow ~^\prime E_{2}^{\ast\ast} $ taking
$x$ to $x$. It follows that $^{\prime}E_{2}^{\ast\ast}$ collapses
proving our claim.\end{proof}

Let us now return to the proof of the theorem; recall that we are
still handling the special case $A=R$. Consider the set of maps
$R(d)\rightarrow R_d$ in the homotopy category of $R$-algebras
over $R$. By Theorem \ref{lp}  this set is an abelian group of
topological derivations of $R(d)$ with values in $\Sigma^{-d}R$.
 Since $R(d)$ is a commutative $R$-algebra there is a canonical splitting
\[THH_R^{\ast}(R(d),R)\simeq R\vee Der_R^{\ast-1}(R(d),R).\] So the
computation of $Der_R^{\ast}(R(d),R)$ reduces to the computation
of topological Hochschild cohomology  $THH^{\ast}_R(R(d),R)$.
Further Corollary \ref{alg} provides an isomorphism
$$THH^{\ast}_R(R(d),R)\cong F_{R(d)}(R,R).$$ It follows from Lemma
\ref{cop} that the spectral sequence
\begin{align*}E^{ij}_2&=Der_{R_\ast}(R(d)_\ast,R_\ast)\\&=
Ext^\ast_{R(d)_\ast}(\Omega_{R(d)_\ast}^{R_\ast},R_\ast)\Rightarrow
{ Der}^\ast_R(R(d), R)\end{align*} collapses and \[{
Der}_R^{\ast-1}(R(d), R)=R_\ast[[x]]/R_\ast.\]

Further, $$[R(d),R]^\ast = R_\ast\oplus \Sigma^d R_\ast =
\Lambda_{R_\ast}(z)$$ where the symbol $z$ has degree $d$ (of
course, we do not claim the existence of any multiplicative
structure). The element $z$ maps the wedge summand $\Sigma^{-d}R$
of the $R$-module $R(d)$ isomorphically to $\Sigma^{-d}R$, and the
other wedge summand maps to zero. It follows from Proposition
\ref{an} that the image of $x\in Der^\ast(R(d),R)$ in
$[R(d),R]^\ast$ under the forgetful map is $z$ (up to an
invertible factor).

In other words we proved that there exists a topological
derivation of $R(d)$ with values in $R$, that is a map in the
homotopy category of $R$-algebras over $R$
\begin{equation}\label {asz} R(d)\rightarrow R_d
\end{equation}
such that the wedge summand $\Sigma^{-d}R$ of $R(d)$ maps
isomorphically onto the corresponding wedge summand of $R_d$.
Therefore the map (\ref{asz}) is a weak equivalence of
$R$-algebras and our theorem is proved (in the special case
$A=R$). To get the general case consider the canonical map
\begin{equation}\label{rty}R(d)\wedge A =F({\bf R}\Sigma^\infty S^d_+,R)\wedge A
\rightarrow F({\bf R}\Sigma^\infty S^d_+,A)=A(d)
\end{equation}
Since $A$ is a $q$-cofibrant $R$-algebra the point-set level smash
product\\ $F({\bf R}\Sigma^\infty S^d_+,R)\wedge A$ represents the
derived smash product. Further (\ref{rty}) is a weak equivalence
since ${\bf R}\Sigma^\infty S^d_+$ is a finite cell $R$-module and
diagram chase shows that this is an $R$-algebra map. So we have
the following equivalences of $R$-algebras: $$A_d\cong R_d\wedge
A\simeq R(d)\wedge A \simeq A(d)$$ With this Theorem \ref{thy} is
proved.\end{proof}

Now suppose that we have another $R$-algebra $B$ and a map
$f:A\rightarrow B$ of $R$-algebras. Then the pair
$(F_{R-alg}(A,B),f)$ is a pointed topological space. This space
turns out to be closely related to $Der_R(A,B)$. We have the
following theorem:
\begin{theorem}\label{osn}For a $q$-cofibrant algebra $A$ and a map of $R$-algebras
$f:A\rightarrow B$ the space of $d$-fold loops
$\Omega^d(F_{R-alg}(A,B),f)$ is weakly equivalent to the space
$\Omega^\infty{\bf Der}_R(A,\Sigma^{-d}B)$.
\end{theorem}
\begin{proof} For two topological spaces $X$ and $Y$ we will denote the
space of maps between them by ${\mathcal T}(X,Y)$ $({\mathcal
T}_\ast(X,Y)$ in the pointed case). Then we have the following
commutative diagram of spaces where both rows are homotopy fibre
sequences: \[\xymatrix{{\mathcal T}_\ast(S^d,F_{R-alg}(A,
B))\ar[r]\ar[d]&{\mathcal T}(S^d,F_{R-alg}(A,
B))\ar[r]\ar[d]&F_{R-alg}(A, B)\ar@{=}[d]\\ ?\ar[r]&
F_{R-alg}(A,B^{S^d})\ar[r]&F_{R-alg}(A,B)}
\]
Here the horizontal rightmost arrows are both induced by the
inclusion of the base point into $S^d$. Since the right and the
middle vertical arrows are weak equivalences (even isomorphisms)
it follows that the map
 ${\mathcal T}_\ast(S^d,Map(A, B))\rightarrow ?$ is a weak equivalence. But
Theorem \ref{thy} tells us that the $R$-algebra  $B^{S^d}$ is
weakly equivalent as an $R$-algebra to $B\vee \Sigma^{-d}B$. In
other words the term $?$ is weakly equivalent to the topological
space of maps $A\rightarrow B\vee \Sigma^{-d}B$ which commute with
the projection onto $A$. Therefore $?$ is weakly equivalent to
$\Omega^\infty{\bf Der}_R(A,\Sigma^{-d}B)$ and our theorem is
proved.\end{proof}
\begin{cor}\label{com}
For a $q$-cofibrant algebra $A$ and a map of $R$-algebras
$f:A\rightarrow B$ there is a bijection between sets
$\pi_d(F_{R-alg}(A,B),f)$ and $Der^{-d}_R(A,B)$ for $d\geq 1$. If
$d\geq 2$ then this bijection is an isomorphism of abelian groups.
\end{cor}
\begin{rem}
One might wonder whether Theorem \ref{osn} remains true in the
context of commutative $S$-algebras. The answer is no. The crucial
point is the weak equivalence of $S$-algebras $S\vee S^{-1}$ and
$S^{S^1}$.  It is clear that $\pi_0S\wedge_{S\vee S^{-1}}S $ is
the divided power ring. However N.Kuhn and M.Mandell proved that
$\pi_0S\wedge_{S^{S^1}}S $ is the ring of numeric polynomials.
Therefore  $S\vee S^{-1}$ and $S^{S^1}$ cannot be weakly
equivalent as commutative $S$-algebras.
\end{rem}
We see, that the space $F_{R-alg}(A, B)$ when looped only once
becomes an infinite loop space. This is somewhat surprising since
$F_{R-alg}(A, B)$ is hardly ever an infinite loop space itself. In
particular the set of connected components of $F_{R-alg}(A, B)$
does not have to be a group, let alone an abelian group. Therefore
the connection between $\pi_0F_{R-alg}(A, B)$ and $Der^0_R(A,B)$
(provided the latter is defined) may be rather weak. For instance
the set of homotopy classes of $A_\infty$ self-maps of the
$p$-completed $K$-theory spectrum is the multiplicative group of
$p$-adic integers whereas the corresponding topological
derivations spectrum  can be proved to be contractible. A
generalization of this example is discussed in author's work
\cite{Laz1}. However there is some evidence for the following
\begin{con} For an $R$-algebra map $f:A\rightarrow B$ the connected component
of $f$ in $F_{R-alg}(A,B)$ is weakly equivalent to the connected
component of $\Omega^\infty{\bf Der}_R(A,B)$. In particular it is
an infinite loop space.
\end{con}
To see why this conjecture has a chance of being true notice that
the Whitehead products in the homotopy groups of $F_{R-alg}(A,B)$
determine via Theorem \ref{osn} various brackets in $Der^*_R(A,B)$
and, for commutative $A$ and $B$ - also in $THH^*_R(A,B)$. No such
brackets have been recorded so far and it seems likely that they
should all vanish. This suggests that the connected component of
$f$ in $F_{R-alg}(A,B)$ is an $H$-space.

There is another interesting question raised by Theorem \ref{osn}.
In recent work \cite{MS} J.McClure and J.Smith introduced the
Gerstenhaber bracket on ${\bf THH}_R(A,A)$. Their work probably
implies the existence of the bracket on ${\bf Der}_R(A,A)$. This
is surely the case if $A$ is commutative since then ${\bf
Der}_R(A,A)$ splits off ${\bf THH}_R(A,A)$ as a wedge summand.
Then via Theorem \ref{osn} a Poisson bracket is defined on
$\pi_*F_{R-alg}(A,A)$ for $*>0$.
\begin{con} The bracket described above agrees with the Whitehead product on
$BF_{R-alg}(A,A)$, the classifying space of the monoid
$F_{R-alg}(A,A)$.
\end{con}

We see that the problem of computing $\pi_0F_{R-alg}(A, B)$
differs sharply from computing higher homotopy groups. This
problem is usually much harder, being essentially nonabelian.
However there is one case when it is possible to give a complete
general answer.

\begin{theorem}\label{eka} Assuming that $R$ is connective  let $A$ be a connective
$q$-cofibrant $R$-algebra, and $B$ a coconnective $R$-algebra
(i.e. $\pi_iB=0$ for $i>0$). Then any $\pi_0R$-algebra map
$\pi_0A\rightarrow \pi_0B$ lifts to a unique $R$-algebra map
$A\rightarrow B$ so that the forgetful map $[A,B]_a\rightarrow
Hom_{\pi_0R-alg}(\pi_0A,\pi_0B)$ is bijective. Moreover the
topological space $F_{R-alg}(A,B)$ is homotopically discrete, i.e.
$\pi_iF_{R-alg}(A,B)=0$ for $i>0$.

Similarly if $A$ and $B$ are both commutative $R$-algebras, where
$B$ is coconnective  and $A$ is $q$-cofibrant and connective then
the forgetful map $[A,B]_c\rightarrow
Hom_{\pi_0R-alg}(\pi_0A,\pi_0B)$ is bijective and the space of
{\it commutative} $R$-algebra maps from $A$ to $B$ is
homotopically discrete.
\end{theorem}
\begin{proof} We will  deal only  with the associative case, the
commutative one being completely analogous. Picking a system of
generators and relations for the $\pi_0R$-algebra $\pi_0A$ we
construct the following pushout diagram in the category of
$R$-algebras:
\begin{equation}\label{df}\begin{array}{ccc}T(\coprod_JR)&\rightarrow&R\\
\downarrow&&\downarrow\\ T(\coprod_IR)&\rightarrow&A^0
\end{array}\end{equation}
Here the sets $I$  and $J$ run respectively through the systems of
generators and relations in $\pi_0A$. There is a canonical
$R$-algebra map from $A^0$ to $A$ that induces an isomorphism on
zeroth homotopy group. The $R$-algebra $A^0$ is the zeroth
skeleton of $A$ in the category of $R$-algebras and (the
$CW$-approximation of) $A$ is obtained from $A^0$ by attaching
$R$-algebra cells in higher dimensions. Then induction up the
$CW$-filtration of $A$ shows that the map $A^0\rightarrow A$
induces a weak equivalence $F_{R-alg}(A,B)\simeq
F_{R-alg}(A^0,B)$.

Further applying the functor $F_{R-alg}(?,B)$ to the diagram
(\ref{df}) we get the following homotopy pullback of topological
spaces: $$\begin{array}{ccc}\prod_J\pi_0B&\leftarrow&pt\\
\uparrow&&\uparrow\\
\prod_I\pi_0B&\leftarrow&F_{R-alg}(A^0,B)\end{array}$$ It follows
that the space $F_{R-alg}(A^0,B)$ is homotopically discrete with
the set of connected components being equal to
$Hom_{\pi_0R-alg}(\pi_0A,\pi_0B)$. \end{proof}

Now let $k$ be an associative ring. Recall that according to
\cite{EKMM}, Proposition IV.3.1 the Eilenberg-MacLane spectrum
$Hk$ admits a structure of an $S$-algebra or a commutative
$S$-algebra if $k$ is commutative. Theorem \ref{eka} shows that
this structure is unique up to a weak equivalence of $S$-algebras
or commutative $S$-algebras. We also have the following evident
corollary which will be used in the next section.
\begin{cor}\label{ek1} Let $A$ be a connective $q$-cofibrant $S$-algebra or
commutative $S$-algebra. Then the topological space of $S$-algebra
maps (or commutative $S$-algebra maps) from $A$ to $H\pi_0A$ is
homotopically discrete and
\[\pi_0F_{S-alg}(A,H\pi_0A)=End_{rings}(\pi_0A).\]\end{cor}

\section{Spaces of multiplicative self-maps of $MU$}
In this section we study the homotopy groups of $A_\infty$-maps
from the complex cobordism spectrum $MU$ into itself. Our main
result here is that any  homotopy multiplicative operation
$MU\rightarrow MU$ lifts canonically to an $S$-algebra map. We
also calculate completely higher homotopy groups of $S$-algebra
maps out of $MU$ into an arbitrary $MU$-algebra $E$. In this
section we work with various homotopy categories and so smash
products and function spectra are understood in the derived sense.

Before we state our main theorem we need to introduce the notion
of $\mathbb{Q}$-commutative $S$-algebras and
 $\mathbb{Q}$-preferred $S$-algebra maps.
\begin{defi}
Let $A$ be an $S$-algebra and denote by $A_\mathbb{Q}$ its
rationalization. We say that $A$ is $\mathbb{Q}$-commutative if
the $A_\mathbb{Q}$ is weakly equivalent as an $S$-algebra to a
{\it commutative} $S$-algebra.
\end{defi}
\begin{rem}
 Later on all $\mathbb{Q}$-commutative $S$-algebras which we encounter will in fact
be commutative. Notice, however, that it is not always the case.
Denote by $S[x_i]$ the free $S$-algebra on the $S$-module
$S_S^{2i}$, the cell approximation of the $2i$-dimensional sphere.
Then clearly $S_\mathbb{Q}[x_i]$ is weakly equivalent to the free
{\it commutative} $S$-algebra on $S_\mathbb{Q}^{2i}$. Therefore
$S[x_i]$ is a $\mathbb{Q}$-commutative $S$-algebra which is not
commutative unless $i=0$.
\end{rem}
Consider two $\mathbb{Q}$-commutative $S$-algebras $A_\mathbb{Q}$
and $B_\mathbb{Q}$. We have the following maps:
\[\xymatrix{k:[A_\mathbb{Q},B_\mathbb{Q}]_{c}\ar[r]&
[A_\mathbb{Q},B_\mathbb{Q}]_{a}&[A,B]_{a}:q\ar[l]}.\] Here $k$ is
the forgetful map and $q$ is induced by rationalization.
\begin{defi} A map $f\in[A,B]_{a}$ is called $\mathbb{Q}$-preferred if
$q(f)\in[A_\mathbb{Q},B_\mathbb{Q}]_{a}$ is in the image of $k$.
Similarly for an $A$-bimodule $M$ which is $\mathbb{Q}$-symmetric
(that is, the square-zero extension $A\vee M$ is
$\mathbb{Q}$-commutative) an $S$-algebra derivation
$d:A\rightarrow A\vee M$ is called $\mathbb{Q}$-preferred if $d$
is $\mathbb{Q}$-preferred as an $S$-algebra map.
\end{defi}
In other words a map of $S$-algebras (or a topological derivation)
is $\mathbb{Q}$-preferred if it lifts to a map (to a derivation)
of commutative $S$-algebras after rationalization.
\begin{theorem}\label{maps} The forgetful map of monoids
\[\xymatrix{[MU,MU]_a\ar[r]& Mult(MU,MU)}\] admits a unique section whose image
consists of $\mathbb{Q}$-preferred $S$-algebra maps.
\end{theorem}
The proof will be given below after a succession of lemmas.
\begin{rem} The set  $Mult(MU,MU)$ is relatively well understood. One can
describe it for example as the set of all $MU_\ast$-algebra maps
\[\xymatrix{MU_\ast MU=MU_\ast[t_1,t_2,\ldots]\ar[r]& MU_\ast}.\]
\end{rem}
Our next result is the computation of topological Hochschild
cohomology of $MU$ with coefficients in an $MU$-algebra $E$. Since
there is a canonical splitting of spectra \[{\bf
THH}_S(MU,E)\simeq E\vee \Sigma^{-1}{\bf Der}_S(MU,E)\] the
combination of this result with Corollary \ref{com} gives a
complete calculation of higher homotopy groups of the based space
$F_{S-alg}(MU,E)$.
\begin{prop}\label{thh}
 For an $MU$-algebra $E$ considered as an $MU$-bimodule the following isomorphism holds
\[THH_S^*(MU,E)\cong\hat{ \Lambda}_{E_\ast}(y_1, y_2,\ldots)\]
where the hat denotes the completed exterior algebra and the
exterior generator $y_i$ has cohomological degree $2i-1$.
\end{prop}
\begin{proof} Consider the topological Hochschild homology $S$-module of
$MU$ with coefficients in $E$, \[{\bf THH}^S(MU,E):=MU\wedge
_{MU\wedge MU}E\cong MU\wedge_{MU\wedge MU}MU\wedge_{MU} E.\] We
have the spectral sequence of $MU_\ast$-algebras
\begin{align*}
E^2_{\ast\ast}&=Tor_{\ast\ast}^{MU_\ast MU}(MU_\ast,MU_\ast)\\& =
 MU_{\ast}\otimes \Lambda(\tilde{y}_1, \tilde{y}_2,\ldots) \Rightarrow
\pi_\ast {\bf THH}^S(MU,MU).
\end{align*}
Since the differentials applied to the exterior generators
$\tilde{y}_i$ are trivial for dimensional reasons we conclude that
it collapses. It follows that \begin{align*}\pi_\ast {\bf
THH}^S(MU,E)&=\pi_\ast {\bf
THH}^S(MU,MU)\otimes_{MU_\ast}E_\ast\\&=E_{\ast}\otimes
\Lambda(\tilde{y}_1, \tilde{y}_2,\ldots).\end{align*} Now the
result for topological Hochschild cohomology follows  by virtue of
the universal coefficients formula and the isomorphism \[{\bf
THH}_S(MU,E) \cong F_E({\bf THH}^S(MU,E),E).\]Proposition
(\ref{thh}) is proved.
\end{proof}
Recall that we are using the notation $S_\mathbb{Q}[x_i]$ for the
free commutative $S$-algebra on the $S$-module
$S_\mathbb{Q}^{2i}$, the rationalized $2i$-sphere $S$-module. The
coefficient ring of $S_\mathbb{Q}[x_i]$ is isomorphic to
$\mathbb{Q}[x_i]$ where the polynomial generator $x_i$ has degree
$2i$. Further denote the infinite smash power
$S_\mathbb{Q}[x_1]^{\wedge\infty}$ by
$S_\mathbb{Q}[x_1,x_2,\ldots]$.
\begin{lem} There is a weak equivalence of commutative $S$-algebras
 \[\xymatrix{S_\mathbb{Q}[x_1,x_2,\ldots]\ar[r]&MU_\mathbb{Q}}.\]
\end{lem}
\begin{proof} The polynomial generators $x_i$ of the ring $
MU_{\mathbb{Q}\ast}=\mathbb{Q}[x_1,x_2,\ldots]$ determine a
collection of maps $S_\mathbb{Q}^{2i}\rightarrow MU_\mathbb{Q}$
and therefore a map of commutative algebras
$S_\mathbb{Q}[x_1]^{\wedge\infty}\rightarrow MU_\mathbb{Q}$ which
is clearly a weak equivalence.\end{proof}
\begin{defi}
Let $E$ be a ring spectrum (in the traditional up to homotopy
sense) with multiplication $m:E\wedge E\rightarrow E$ and $M$ an
$E$-bimodule spectrum with the left action $m_l:E\wedge
M\rightarrow M$ and the right action $m_r:M\wedge E\rightarrow M$.
We say that a map $f:E\rightarrow M$ is a primitive operation if
$f\circ m$ and $m_r\circ (f\wedge id)+ m_l\circ (id \wedge f)$ are
homotopic as maps from $E\wedge E$ to $M$. The set of all
primitive operation from $E$ to $M$ is denoted by $Prim(E,M)$
\end{defi}
\begin{rem} Perhaps it is more natural to use the term `derivation' instead of
`primitive operation' but this term is already overworked in this
paper.
\end{rem}
The next lemma provides a description of topological derivations
of $MU$ with coefficients in $H\mathbb{Z}$, the integral
Eilenberg-MacLane spectrum.
\begin{lem}\label{ghj} There is the following isomorphism of graded abelian groups:
\[ Der_S^\ast(MU, H\mathbb{Z})\cong
\Lambda^{\ast-1}(y_1,y_2\ldots)/\mathbb{Z}.\] Under the forgetful
map \[\xymatrix{l: Der_S^\ast(MU, H\mathbb{Z})\ar[r]& [MU,
H\mathbb{Z}]^\ast=Hom(\mathbb{Z}[t_1,t_2,\ldots],\mathbb{Z})}\]
the elements $y_i\in Der_S^{2i-2}(MU, H\mathbb{Z})$ correspond to
the derivations $\partial_{t_i}$ evaluated at 0. Moreover the
elements $y_i$ are $\mathbb{Q}$-preferred topological derivations.
\end{lem}
\begin{proof} We have the spectral sequence \begin{align*}
Der^\ast(H\mathbb{Z}_\ast MU,
\mathbb{Z})&=\Lambda^{\ast-1}(y_1,y_2\ldots)/\mathbb{Z}\\&\Rightarrow
Der_{H\mathbb{Z}}^\ast(H\mathbb{Z}\wedge MU,
H\mathbb{Z})=Der_S^\ast(MU, H\mathbb{Z}).\end{align*} This
spectral sequence clearly collapses. Next using Proposition
\ref{an} we see that the image of the element $y_i$ under the
forgetful map $l$ in the group \[[H\mathbb{Z}\wedge MU,
H\mathbb{Z}]^*_{H\mathbb{Z}-mod}=[MU,H\mathbb{Z}]^*=Hom^*(\mathbb{Z}[t_1,t_2,\ldots],
\mathbb{Z})\] is precisely the algebraic derivation
$\partial_{t_i}$ evaluated at $0$ (up to elements of higher
filtration). Since this image is contained in the subgroup of
primitive operations $MU\rightarrow H\mathbb{Z}$ none of these
elements of higher filtration are present.

To see that $y_i$ are $\mathbb{Q}$-preferred derivations let us
introduce the notation $CDer^\ast(MU_\mathbb{Q},H\mathbb{Q})$ to
denote topological derivations of $MU_\mathbb{Q}$ with values in
$H\mathbb{Q}\simeq S_\mathbb{Q}$ in the category of commutative
$S$-algebras. (These commutative derivations are also known as
topological Andre-Quillen cohomology, cf.\cite{Bas}).
 Then since $MU_\mathbb{Q}$ is
a free commutative $S$-algebra we see immediately that
$$CDer^\ast(MU_\mathbb{Q},H\mathbb{Q})=Der^\ast(MU_{\mathbb{Q}\ast},\mathbb{Q})=
\mathbb{Q}<\bar{\partial}_{x_1},\bar{\partial}_{x_2},\ldots>$$ the
right hand side being the set of derivations (in the usual
algebraic sense) of the algebra $MU_{\mathbb{Q}\ast}$ with values
in the rational numbers. Here we denoted by $\bar{\partial}_{x_i}$
the standard derivation $\partial_{x_i}$ of the ring
$MU_\mathbb{Q}=\mathbb{Q}[x_1,x_2,\ldots]$ composed with
evaluation at zero.

On the other hand $Der_S^\ast(MU_\mathbb{Q},H\mathbb{Q})\cong
\Lambda_\mathbb{Q}^{\ast-1}(y_1,y_2\ldots)/\mathbb{Q}$. We need to
prove therefore that the forgetful map
\[\xymatrix{CDer^\ast(MU_\mathbb{Q},H\mathbb{Q})\ar[r]&
Der_S^\ast(MU_\mathbb{Q},H\mathbb{Q})}\] sends the elements
$\bar{\partial}_{x_i}$ to $y_i$.

Since  the commutative $S$-algebra $H\mathbb{Q}[x_i]\simeq
S_\mathbb{Q}[x_i]$
 is free as a commutative $S$-algebra as well as an (associative) $S$-algebra it follows
that
$$CDer^\ast(H_\mathbb{Q}[x_i],H\mathbb{Q})=Der^\ast(\mathbb{Q}[x_i],\mathbb{Q})=\mathbb{Q}<\bar{\partial}_{x_i}>.$$
There is a unique map of commutative $S$-algebras
$MU_\mathbb{Q}\rightarrow H\mathbb{Q}[x_i]$ which corresponds to
quotienting out the ideal
$(x_1,x_2,\ldots,x_{i-1},\hat{x}_i,x_{i+1},\ldots)$ in the
coefficient ring of $MU_\mathbb{Q}$. We have the following
commutative diagram:
\[\xymatrix{CDer^\ast(H\mathbb{Q}[x_i],H\mathbb{Q})\ar^\cong[r]\ar[d]&
Der^\ast(\mathbb{Q}[x_i],\mathbb{Q}))\ar[d]\\
CDer^\ast(MU_\mathbb{Q},H\mathbb{Q})\ar[r]
&Der_S^\ast(MU_\mathbb{Q},H\mathbb{Q})}\] from which it is clear
that the image of $\bar{\partial}_{x_i}$ in
$Der^\ast(MU_\mathbb{Q},H\mathbb{Q})$ is $y_i$ and the lemma is
proved.\end{proof}
\begin{cor} \label{set}The set of $\mathbb{Q}$-preferred derivations of
$MU$ with values in $H\mathbb{Z}$ maps bijectively onto the set of
primitive cohomology operations $MU\rightarrow H\mathbb{Z}$ under
the forgetful map $Der_S^\ast(MU,H\mathbb{Z})\rightarrow
[MU,H\mathbb{Z}]^\ast$.
\end{cor}
\begin{lem}\label{cot} Let $A$ be an $S$-algebra and $B$ a commutative $S$-algebra.
Suppose that $B$ has a structure of an $A$-bimodule via a map of
$S$-algebras $A\rightarrow B$. Then ${\bf THH}_S(A,B)$ has a
structure of a $B$-bimodule and there is a canonical splitting of
$B$-bimodules ${\bf THH}_S(A,B)\simeq B\vee \Sigma^{-1}{\bf
Der}_S(A,B)$.
\end{lem}
\begin{proof} Consider the following sequence of $S$-algebra maps:
\[B\longrightarrow {\bf THH}_S(B,B)\cong F_{B\wedge
B^{op}}(B,B)\longrightarrow
 F_{B\wedge A^{op}}(B,B)\cong {\bf THH}_S(A,B).\]
The first map exists because $B$ is commutative, the middle map is
induced by the $S$-algebra map $A\rightarrow B$ and the last
equivalence is Corollary \ref{alg}. The composite map
$B\rightarrow{\bf THH}_S(A,B)$ supplies ${\bf THH}_S(A,B)$ with a
structure of a $B$-bimodule and splits the canonical map ${\bf
THH}_S(A,B)\rightarrow B$. \end{proof} Let us introduce the
notation $MU_n$ for the $n$th Postnikov stage of $MU$. Then $MU_n$
is an $S$-algebra (even a commutative $S$-algebra).
\begin{lem}\label{uop} There is the following weak equivalence of $ MU_{n\ast}$-modules:
\[THH_S^\ast(MU,MU_n)\cong \Lambda_{MU_{n\ast}}(y_1,y_2,\ldots)\]
where $MU_n$ is considered as an $MU$-bimodule via any (not
necessarily central) map of $S$-algebras $f:MU\rightarrow MU_n$.
\end{lem}
\begin{proof} Since $MU_n$ is a commutative $S$-algebra the multiplication
map \[\xymatrix{MU_n\wedge MU_n\ar^m[r]&MU_n}\] is an $S$-algebra
map. Therefore the composition \[\xymatrix{MU\wedge
MU_n\ar^{f\wedge id}[r]& MU_n\wedge MU_n \ar^-m[r]&MU_n}\] is also
an $S$-algebra map. This gives  the following weak equivalence of
$S$-modules: \begin{align*}{\bf THH}_S(MU,MU_n)&\cong F_{MU\wedge
MU_n}(MU_n,MU_n)\\&\cong F_{MU_n}(MU_n\wedge_{MU\wedge
MU_n}MU_n,MU_n).\end{align*} Therefore it is enough to show that
$$\pi_\ast MU_n\wedge_{MU\wedge MU_n}MU_n=\Lambda_{MU_{n\ast}}
(\tilde{y}_1,\tilde{y}_2\ldots).$$ (The exterior generators
$y_i\in THH_S^\ast(MU,MU_n) $ will be dual to $\tilde{y}_i$).

 Consider the spectral sequence
\begin{equation}\label{spe}Tor_{\ast\ast}^{MU_\ast MU_n}(MU_{n\ast},MU_{n\ast})=
\Lambda_{MU_{n\ast}}(\tilde{y}_1,\tilde{y}_2,\ldots)\Rightarrow
\pi_\ast MU_n\wedge_{MU\wedge MU_n}MU_n \end{equation} This
spectral sequence is {\it not} multiplicative since the map
$f:MU\rightarrow MU_n$ may not be central. However it is a
spectral sequence of $MU_{n\ast}$-modules.

Let us introduce  another spectral sequence
\begin{equation}\label{sp1}Tor_{\ast\ast}^{MU_\ast MU}(MU_{\ast},MU_{\ast})=
\Lambda_{MU_{\ast}}(\tilde{y}_1,\tilde{y}_2,\ldots)\Rightarrow
\pi_\ast MU\wedge_{MU\wedge MU}MU \end{equation} Then the map
$f:MU\rightarrow MU_n$ induces a map of spectral sequences
(\ref{sp1})$\rightarrow$(\ref{spe}). Further the spectral sequence
(\ref{sp1}) is multiplicative and collapses for that reason.
Therefore in (\ref{spe}) all elements of the form
$\tilde{y_{i_1}}\wedge\tilde{y_{i_2}}\wedge\ldots \tilde{y_{i_k}}$
are permanent cycles and it follows that (\ref{spe}) collapses.
Lemma \ref{uop} is proved.\end{proof}

Suppose as before that we have an $MU$-bimodule structure on
$MU_{n+1}$ via some $S$-algebra map $f:MU\rightarrow MU_{n+1}$.
Composing $f$ with the canonical map in the Postnikov tower
$p_n:MU_{n+1}\rightarrow MU_n$ we get an $MU$-bimodule structure
on $MU_{n}$ also. Then we have the following
\begin{cor} \label{pol}The induced map
$Der_S^\ast(MU,MU_{n+1})\longrightarrow Der_S^\ast(MU,MU_{n})$ is
onto.
\end{cor}
\begin{proof} Indeed, the map
\begin{align*}THH_S^{\ast}(MU,MU_{n+1})&=\Lambda_{MU_{n+1\ast}}(y_1,y_2,\ldots)\\
&\rightarrow
\Lambda_{MU_{n\ast}}(y_1,y_2,\ldots)=THH_S^{\ast}(MU,MU_{n})\end{align*}
is clearly onto and our claim follows from Lemma
\ref{cot}.\end{proof}

{\it Proof of Theorem \ref{maps}.} We start by outlining the
strategy of the proof. Take a multiplicative operation $f\in
Mult(MU,MU)$.
 Define  $f_n:MU\rightarrow MU_n$ as the composition
\[\xymatrix{MU\ar^f[r]& MU\ar^{p_n}[r]&MU_n}.\]
(Recall that we denoted by  $p_n$ the canonical projection onto
the $n$th Postnikov stage.)
 Then by Corollary \ref{ek1} the map
$f_0:MU\rightarrow MU_0=H\mathbb{Z}$ is homotopic to a unique
$S$-algebra map which we will denote by
 $\tilde{f}_0$. Proceeding by induction assume
that there exists a unique $\mathbb{Q}$-preferred $S$-algebra map
$\tilde{f}_n:MU\rightarrow MU_n$ which is homotopic to $f_n$ when
considered as a map of $S$-modules. We will see that
\begin{itemize}\item $\tilde{f}_n$ admits a $\mathbb{Q}$-preferred
lifting to an $S$-algebra map $MU\rightarrow MU_{n+1}$ and \item
there is a one-to-one correspondence between such liftings and the
set of liftings of $ \tilde{f}_n$ to a \emph{homotopy}
multiplicative map $MU\rightarrow MU_{n+1}$.\end{itemize} In
particular $f_{n+1}$ being one of such liftings can be realized as
a $\mathbb{Q}$-preferred $S$-algebra lifting in a unique fashion.

We now proceed to realize the above program in detail. The first
thing is to show that there exists a lifting of $\tilde{f}_n$ in
the category of $S$-algebras. The homotopy fibre sequence
\begin{equation}\label{gos}
\xymatrix{\Sigma^{n+1}H\pi_{n+1}MU\ar[r]& MU_{n+1}\ar[r]&
MU_n}\end{equation} is a topological singular extension by Theorem
\ref{pos}. Then Theorem {\ref{ext}} tells us that the obstruction
to an $S$-algebra lifting $\tilde{f}_n$ to $MU_{n+1}$ is a certain
element  $\sigma\in Der^0_S(MU, \Sigma H\pi_{n+1}MU)$. More
precisely, the extension (\ref{gos}) is associated with a
derivation
\[\xymatrix{d:MU_n\ar[r]& MU_n\vee \Sigma^{n+2}H\pi_{n+1}MU}\] and
$\sigma:MU\rightarrow MU_n \vee \Sigma^{n+2}H\pi_{n+1}MU$ is the
composition of $d$ with $p_n:MU\rightarrow MU_n$.

Furthermore, notice that the set of $S$-algebra maps
$MU_\mathbb{Q}\rightarrow MU_{\mathbb{Q}n+1}$ is in bijective
correspondence with the set of ring maps
$MU_{\mathbb{Q}\ast}\rightarrow MU_{\mathbb{Q}n+1\ast}$.
 Since $MU_{\mathbb{Q}\ast}$ is a polynomial algebra we see that a lift of
 $\tilde{f_n}$ does exists after
rationalization. Therefore the image of $\sigma$ in $Der^0_S(MU,
\Sigma^{n+2} H\pi_{n+1}MU_\mathbb{Q})$ is zero. But
~$H\pi_{n+1}MU$ is a wedge of suspensions of $H\mathbb{Z}$ and
according to Lemma \ref{ghj} the abelian group $Der^0_S(MU,
H\pi_{n+1}MU)$ has no torsion. It follows that $\sigma=0$ as an
element in the group $Der^0_S(MU, \Sigma H\pi_{n+1}MU)$ and a lift
of $\tilde{f}_n$ exists integrally (though not necessarily
$\mathbb{Q}$-preferred). By Theorem \ref{ext} the homotopy fibre
of the map \[\xymatrix{F_{S-alg}(MU,MU_{n+1})\ar[r]&
F_{S-alg}(MU,MU_{n})}\] taken over ~ the ~ point ~~
$\tilde{f}_n\in F_{S-alg}(MU,MU_{n}) $  is ~ weakly ~ equivalent~
to~ the~ zeroth space of the spectrum ${\bf Der}_S(MU,
\Sigma^{n+1}H\pi_{n+1}MU)$.

Therefore denoting by $[MU, MU_{n+1}]_a^{lift}\subset [MU,
MU_{n+1}]_a$ the set of homotopy classes of $S$-algebra maps
$MU\rightarrow MU_{n+1}$ lifting $\tilde{f}_n$ we have the
following long exact sequence:
\begin{align*}
\ldots\longrightarrow Der_S^{-1}(MU,\Sigma^{n+1}H\pi_{n+1}MU)
&\longrightarrow  \pi_1 F_{S-alg}(MU, MU_{n+1})
\\ \\
 \longrightarrow\pi_1 F_{S-alg}(MU, MU_{n})  &\longrightarrow
 Der_S^{0}(MU,\Sigma^{n+1}H\pi_{n+1}MU)
 \\ \\\longrightarrow [MU,
MU_{n+1}]^{lift}_a&\longrightarrow pt \end{align*} which is the
same (by Theorem \ref{osn}) as the long exact sequence
\begin{align*}\ldots\longrightarrow
Der^{-1}_S(MU,\Sigma^{n+1}H\pi_{n+1}MU)&\longrightarrow
Der_S^{-1}(MU, MU_{n+1})\\ \\ \longrightarrow
Der^{-1}_S(MU,MU_{n})&\longrightarrow
Der^{0}_S(MU,\Sigma^{n+1}H\pi_{n+1}MU)\\ \\ \longrightarrow [MU,
MU_{n+1}]_a^{lift}&\longrightarrow pt.
\end{align*}
By Corollary \ref{pol} the map
\[\xymatrix{Der^{-1}_S(MU, MU_{n+1})\ar[r]&Der_S^{-1}(MU, MU_{n})}\] is
onto and we conclude that the map
\[\xymatrix{Der^{0}(MU,\Sigma^{n+1}H\pi_{n+1}MU)\ar[r]& [MU,
MU_{n+1}]_a^{lift}}\] is bijective. So the set of all possible
lifts of $\tilde{f}_k$ is in one-to-one correspondence with
elements in the group $Der^{0}_S(MU,\Sigma^{n+1}H\pi_{n+1}MU)$.
 Clearly the set of all $\mathbb{Q}$-preferred lifts corresponds under this isomorphism to
the set of $\mathbb{Q}$-preferred topological derivations of $MU$
with values in $\Sigma^{n+1}H\pi_{n+1}MU$. By Corollary \ref{set}
these $\mathbb{Q}$-preferred derivations are identified with the
set of primitive operations from $MU$ to
$\Sigma^{n+1}H\pi_{n+1}MU$. So we established a one-to-one
correspondence between the set of $\mathbb{Q}$-preferred lifts of
$\tilde{f_n}$ and $Prim(MU,\Sigma^{n+1}H\pi_{n+1}MU)$.

Now we examine the question of lifting the map $\tilde{f}_n$  up
to homotopy to a homotopy multiplicative map $MU\rightarrow
MU_{n+1}$. Clearly the homotopy class of any map of $S$-modules
$MU\rightarrow MU_k$ is determined by its rationalization, i.e,
the rationalization map $[MU,MU_k]\rightarrow
[MU_\mathbb{Q},MU_{\mathbb{Q}k}]$ is injective. It follows that
the map \[\xymatrix{Mult(MU,MU_k)\ar[r]&
Mult(MU_\mathbb{Q},MU_{\mathbb{Q}k})}\] is also injective.

Further we have the following bijection (for any $k$)
$$Mult(MU_\mathbb{Q},MU_{k\mathbb{Q}})\cong
Hom_{rings}(MU_{\mathbb{Q}\ast},MU_{\mathbb{Q}k\ast})$$ Therefore
there is a short exact sequence
\begin{align}\label{gt}\nonumber 0&\longrightarrow
Prim(MU_\mathbb{Q},\Sigma^{n+1}H\pi_{n+1}
MU_{\mathbb{Q}})=Der(MU_{\mathbb{Q}\ast},\pi_{n+1}
MU_{\mathbb{Q}}) \\ \nonumber \\&\longrightarrow
Mult(MU_\mathbb{Q},MU_{\mathbb{Q}n+1})\longrightarrow
Mult(MU_\mathbb{Q},MU_{\mathbb{Q}n})\longrightarrow pt
\end{align}
Of course the last three terms are only sets. The exactness here
means that $Mult(MU_\mathbb{Q},MU_{\mathbb{Q}n+1})$ has a faithful
action of the group $Prim(MU_\mathbb{Q},\Sigma^{n+1}H\pi_{n+1}
MU_{\mathbb{Q}})$ so that the quotient is isomorphic to
$Mult(MU_\mathbb{Q},MU_{\mathbb{Q}n})$.

Consider the diagram of fibre sequences
\[\xymatrix{\Sigma^{n+1}H\pi_{n+1}
MU\ar[r]\ar[d]&MU_{n+1}\ar[r]\ar[d]& MU_n\ar[d]\\
\Sigma^{n+1}H\pi_{n+1}
MU_{\mathbb{Q}n+1}\ar[r]&MU_{\mathbb{Q}n+1}\ar[r]
&MU_{\mathbb{Q}n}}\] Taking into account the fact that
$MU_{k+1}^\ast MU$ surjects onto $MU_{k}^\ast MU$ for any $k$ we
obtain a map of short exact sequences \[\xymatrix{0 \rightarrow
[MU,\Sigma^{n+1}H\pi_{n+1} MU]\ar[r]\ar[d]
&[MU,MU_{n+1}]\ar[r]\ar[d]& [MU,MU_n]\ar[d]\rightarrow 0\\
0\rightarrow [MU_\mathbb{Q},\Sigma^{n+1}H\pi_{n+1}
MU_{\mathbb{Q}}]\ar[r]&[MU_\mathbb{Q},MU_{\mathbb{Q}n+1}]\ar[r]&
[MU_\mathbb{Q},MU_{\mathbb{Q}n}]\rightarrow 0}\] Notice that all
downward maps are injections. Combining this with (\ref{gt}) we
find that there is a short exact sequence \begin{align*}
0&\longrightarrow Prim(MU,\Sigma^{n+1}H\pi_{n+1} MU)\\
&\longrightarrow Mult(MU,MU_{n+1})\longrightarrow
Mult(MU,MU_{n})\longrightarrow 0\end{align*} That is the
indeterminacy in lifting the map $\tilde{f}_n:MU\rightarrow MU_n$
is precisely the set of primitive cohomology operations
$Prim(MU,\Sigma^{n+1}H\pi_{n+1} MU)$. We see that the set of lifts
of the map $\tilde{f}_n$ to a a homotopy multiplicative map is in
one-to-one correspondence with $\mathbb{Q}$-preferred $S$-algebra
lifts. This completes the inductive step and shows that the
original homotopy multiplicative map $f:MU\rightarrow MU$ can be
improved in a unique way to a $\mathbb{Q}$-preferred $S$-algebra
map.

So we succeded in finding a section $i:Mult(MU,MU)\rightarrow
[MU,MU]_a$ of the forgetful map $j:[MU,MU]_a\rightarrow
Mult(MU,MU)$ so that the image of $i$ consists of
$\mathbb{Q}$-preferred $S$-algebra self-maps of $MU$. To see that
$i$ respects composition notice that for $f,g\in Mult(MU,MU)$ the
$S$-algebra map $i(f)\circ i(g)$ is $\mathbb{Q}$-preferred and
$j(i(f)\circ i(g))=f\circ g$. Since there is a unique
$\mathbb{Q}$-preferred $S$-algebra self-map whose image under $j$
is $f\circ g$ we conclude that $i(f)\circ i(g)=i(f\circ g)$. With
this the proof of
 Theorem \ref{maps}
is completed.
\begin{rem}
Using the Bousfield-Kan mapping space spectral sequence (cf.
\cite{Bou}) it is possible to calculate the set of all $S$-algebra
self-maps of $MU$. However this approach leads to the
identification of $[MU,MU]_a$ only as a set, not as a monoid. It
seems that the monoid structure on $[MU,MU]_a$ should be related
to the Gerstenhaber bracket on $THH^*_S(MU,MU)$
\end{rem}
Now consider an $S$-algebra $E$ with a fixed map of $S$-algebras
$f:MU\rightarrow E$ Suppose that $E$  satisfies the following
condition:

$(S)$ The unit map $f:MU_\ast\rightarrow E_\ast$ is surjective.
\begin{rem}
In \cite{Goe} it was proved that a rather broad class of
$C$-oriented spectra (namely those which are obtained by killing
any regular ideal in the ring $MU$) can be supplied with
$MU$-algebra structures. For this class of spectra the condition
$(S)$ is obviously satisfied.
\end{rem}
\begin{cor}
For an $S$-algebra $E$ satisfying the condition $(S)$ any
multiplicative operation $MU\rightarrow E$ can be lifted (perhaps
in a non-unique way) to an $S$-algebra map.
\end{cor}
\begin{proof} The condition $(S)$ guarantees that the map \[
Mult(MU,MU)\rightarrow Mult(MU,E)\] induced by the given map
$f:MU\rightarrow E$ is surjective. In other words any
multiplicative operation $g:MU\rightarrow E$ can be represented as
a composition $h\circ f$ where $h\in Mult(MU,MU)$. Since $h$ can
be lifted to an $S$-algebra self-map of $MU$ our claim
follows.\end{proof}

As another consequence of Theorem \ref{maps} we will show that the
$p$-local Brown-Peterson spectrum $BP$ is an $A_\infty$-retract of
$MU_{(p)}$, the spectrum $MU$ localized at $p$. Recall from e.g.
\cite{rav}, 4.1 that there exists a multiplicative cohomology
operation $g:MU_{(p)}\rightarrow MU_{(p)}$ which is idempotent and
whose image is the $p$-local spectrum $BP$.

\begin{theorem}\label{BP} There exists an $S$-algebra map $f:MU_{(p)}\rightarrow BP$ which has
a right inverse $S$-algebra map $h:BP\rightarrow MU_{(p)}$.
\end{theorem}
\begin{proof} According to Theorem \ref{maps} the multiplicative operation
$g:MU_{(p)}\rightarrow MU_{(p)}$ determines a map of $S$-algebras
which we will denote by the same letter. Without loss of
generality we can assume $g$ to be a $q$-cofibration of
$S$-algebras. Consider the diagram in the  category of
$S$-algebras: \[\xymatrix{
MU_{(p)}\ar^{g}[r]\ar^g[d]&MU_{(p)}\ar^{g}[r]\ar^g[d]&\ldots\\
MU_{(p)}\ar^{id}[r]&MU_{(p)}\ar^{id}[r]&\ldots}\] Each square in
this diagram is commutative since the operation $g$ is idempotent.
The colimit of the upper row taken in the category of $S$-algebras
coincides with the colimit taken in the category of spectra by
Cofibration Hypothesis (\cite{EKMM},VII.4) and both are equivalent
to $BP$. (That shows that $BP$ is an $S$-algebra). Now the map
$f:MU_{(p)}\rightarrow BP$ is just the canonical map to the
colimit. Next the colimit of the lower row is obviously $MU_{(p)}$
and therefore there exists an $S$-algebra map $h:BP\rightarrow
MU_{(p)}$. It follows that $f\circ h:BP\rightarrow BP$ is
homotopic to the identity and Theorem \ref{BP} is
proved.\end{proof}
\begin{rem} It can be shown (cf. \cite{Laz}, \cite {BaJ}, \cite{Goe})
that $BP$ actually supports a structure of an $MU$-algebra.
\end{rem}
\begin{rem}It seems natural to conjecture that any $\mathbb{Q}$-preferred $S$-algebra
self-map of $MU$ lifts to a {\it commutative} $S$-algebra
self-map. This conjecture, if true, would imply the existence of a
canonical $E_\infty$ ring structure on $BP$, a long-standing
problem posed by P.May. The first (to author's knowledge) serious
attack on this problem was undertaken by I. Kriz in his 1993
preprint \cite {Kri}. This paper  inspired much activity in the
area, however it is still regarded as a program for further work
rather than a definitive solution.
\end{rem}
  Even though we don't know whether $BP$ is a commutative $S$-algebra
we can use Theorem \ref{BP} to compute homotopy classes of
$A_\infty$-maps out of $BP$.
\begin{cor} \label{maps1}
For an $S$-algebra $E$ satisfying the condition  $(S)$ every
multiplicative operation $BP\rightarrow E$
 lifts to an $A_\infty$ ring
map $BP\rightarrow E$ (perhaps in a non-unique way).
\end{cor}
\begin{proof} The composition of the multiplicative operation
$BP\rightarrow E$ with the canonical projection $MU\rightarrow BP$
determines a multiplicative operation $MU\rightarrow E$. This
operation lifts to an $S$-algebra map. Composing this $S$-algebra
map with the splitting map $BP\rightarrow MU$ (which we know is an
$S$-algebra map by Theorem \ref{BP}) we find the desired
$S$-algebra map $BP\rightarrow E$.\end{proof}

$Acknowledgement.$ The author would like to thank P.Goerss and
N.Kuhn for many stimulating discussions.

\end{document}